\documentclass[11pt]{amsart}

\usepackage{fullpage, amsmath,amssymb}
\usepackage{url}
\usepackage{enumerate}
\usepackage{graphicx}
\usepackage{color}
\usepackage{tikz}
\usepackage{tikz-3dplot}
\usepackage[labelformat=empty]{subfig}
\usepackage[euler-digits]{eulervm}
\usepackage[all]{xy}
\usepackage{tikz}
\usepackage{tikzcd}
\usetikzlibrary{matrix,arrows}
\usepackage{hyperref}

\author[D.Karp]{Dagan Karp}
\address{Dagan Karp \\ Department of Mathematics, Harvey Mudd College, 340 E Foothill Boulevard, Claremont, CA 91711, USA}
\email{dagan.karp@hmc.edu}
\author[D.Ranganathan]{Dhruv Ranganathan}
\address{Dhruv Ranganathan \\ Department of Mathematics, Yale University, 10 Hillhouse Avenue, New Haven, CT 06511, USA}
\email{dhruv.ranganathan@yale.edu}
\date{\today}

\newtheorem{theorem}{Theorem}

\newtheorem{lemma}[theorem]{Lemma}

\newtheorem{definition}[theorem]{Definition}

\newtheorem{example}[theorem]{Example}
\newtheorem{remark}[theorem]{Remark}

\newcommand{\PP}{\mathbb{P}}

\newcommand{\ZZ} {{\mathbb Z}}

\newcommand{\Cross}{\!  \times \! }

\newcommand{\M}{\overline{{M}}}

\newcommand{\ta}{\tilde{a}}

\newcommand{\td}{\tilde{d}}

\begin{document}

\title{Gromov-Witten Theory of $\PP^{1}\Cross \PP^{1}\Cross \PP^{1}$}

\pagestyle{plain}

\begin{abstract}
We use elementary geometric techniques to exhibit an explicit equivalence between certain sectors of the Gromov-Witten theories of blowups of $\PP^{1}\Cross \PP^{1}\Cross \PP^{1}$ and $\PP^{3}$.   
In particular, we prove that the all genus, virtual dimension zero
Gromov-Witten theory of the blowup of $\PP^3$ at points
coincides with that of the blowup at points of $\PP^1 \Cross
\PP^1\Cross\PP^1$,  for non-exceptional classes. 
We observe a toric symmetry of the Gromov-Witten theory of
$\PP^1\Cross\PP^1\Cross\PP^1$ analogous and intimately related to
Cremona symmetry of $\PP^3$. Enumerative applications are given.
\end{abstract}

\maketitle

\section{Introduction}\label{sec: introduction}

In this note, we describe a relationship between the Gromov--Witten
theories of the blowups at points of $\PP^3$ and $(\PP^1)^{3}$. This
provides a precise correspondence between curve counts on these
spaces, in the cases considered.  We prove this relationship by
appealing to the toric geometry of the permutohedron.  

\subsection{Results} Let $X$ be the blowup of $\PP^3$ at $6$ points
and $\tilde X$ the blowup of $(\PP^1)^3$ at $4$ points. Let $h\in
H_2(X;\ZZ)$ be the pullback of the class of a line from $\PP^3$, and let
$e_i$ be the class of a general line in the $i^{\textrm{th}}$ exceptional
divisor in $X$. Similarly, let $\tilde h_i$ be the pullback to $\tilde
X$ of the class of the $i^{\textrm{th}}$ line in $(\PP^1)^{3}$ and
$\tilde e_j$ the class of a general line in the  $j^{\textrm{th}}$
exceptional divisor in $\tilde X$. Qualitatively, our first result
states that the all genus virtual dimension zero Gromov--Witten
invariants of $X$ and $\tilde X$ coincide for classes that are far
from the exceptional locus.  

\begin{theorem}\label{thm: p3p1-theorem}
Let $X$ and $\tilde X$ be as above. If $\beta = dh - a_1 e_1 - \ldots
- a_6e_6 \in A_1(X)$ with $(a_5,a_6) \neq (0,0)$, then for any genus $g$, we have
\[
\langle \ \rangle^X_{g,\beta} = \langle \ \rangle^{\tilde X}_{g,\tilde \beta},
\]
where $\tilde \beta = \sum_{i=1}^{3} \tilde d_{i} \tilde h_{i} -\sum_{i=1}^{4} \tilde a_i \tilde e_i$ and the coefficients of $\beta$ and $\tilde \beta$ are related by
\begin{eqnarray*}
\tilde d_{1} &=& d-a_2-a_3 \\
\tilde d_{2} &=& d-a_1-a_3 \\
\tilde d_{3} &=& d-a_1-a_2 \\
\tilde a_1 &=& a_4 \\
\tilde a_2 &=& d-a_1-a_2-a_3 \\
\tilde a_3 &=& a_{5} \\
\tilde a_4 &=& a_{6} .
\end{eqnarray*}
\end{theorem}

The Gromov--Witten theory of $X$ admits a symmetry induced by the Cremona transform (see Section~\ref{sec: further-discussion} below). We find that the Gromov--Witten invariants of $\tilde X$ also admit a symmetry, described as follows.

\begin{theorem}\label{thm: p1-symmetry-theorem}
Let $\tilde X$ be the blowup of $(\PP^1)^3$ at $4$ points. If
$\tilde \beta = \sum_{1\leq j \leq 3} \tilde d_{j} \tilde h_{j}
-\sum_{i=1}^4 \tilde a_i \tilde e_i$, and $\{a_3,a_4\}\neq \{0\}$, we
have 
\[
\langle \ \rangle^{\tilde X}_{g,\tilde \beta} = \langle \ \rangle^{\tilde X}_{g,\tilde \beta'},
\]
where $\tilde{\beta}' = \sum_{1\leq i<j \leq 3} \tilde d'_{j} \tilde h_{j} -\sum_{i=1}^4 \tilde a'_i \tilde e_i$ has coefficients given by 
\begin{eqnarray*}
\tilde d_{1}' &=& \tilde d_{1}+ \tilde d_{3}-\tilde a_1-\tilde a_2\\
\tilde d_{2}' &=& \tilde d_{2}+ \tilde d_{3}-\tilde a_1-\tilde a_2\\
\tilde d_{3}' &=& \tilde d_{3}\\
\tilde a_1'&=& \tilde d_{3}-\tilde a_2\\
\tilde a_2' &=& \tilde d_{3}-\tilde a_1\\
\tilde a_3' &=& \tilde a_4 \\
\tilde a_4' &=& \tilde a_3. 
\end{eqnarray*}
\end{theorem}

\subsubsection{Enumerative consequences} One may relate invariants of $X$ and $\tilde{X}$ to invariants of $\PP^{3}$ and $(\PP^{1})^{3}$ using the following result of Bryan-Leung~\cite{BL}, generalizing a result of Gathmann~\cite{Ga01}.
Let $Y$ be a smooth variety and $\pi :\hat Y\to Y$ the blowup of $Y$ at a point. Let $\beta\in A_1(Y)$ and $\hat \beta = \pi^!(\beta)$. Then,
\[
\langle p \rangle^Y_{g,\beta} = \langle \ \rangle^{\hat Y}_{g,\hat\beta-\hat e},
\]
where $\hat e$ is the class of a line in the exceptional locus, and $\pi^!(\beta) = [\pi^\star[\beta]^{\textnormal{PD}}]^{\textnormal{PD}}$.

\begin{example}\label{ex: enumerative}
How many rational curves in $(\PP^{1})^{3}$ of class $h_{1}+h_{2}+h_{3}$ pass through three general points? We compute 
\begin{align*}
\langle p^{3} \rangle_{0,h_{1}+h_{2}+h_{3}}^{(\PP^{1})^{3}} &= 
 \langle \; \rangle_{0,h_{1}+h_{2}+h_{3}-e_{1}-e_{3}-e_{4}}^{\tilde{X}} \\
&= \langle \; \rangle_{0,3h -e_{1}-\dotsb -e_{6}}^{X}\\
&= \langle p^{6} \rangle_{0,3h}^{\PP^{3}}\\
&=1.
\end{align*}
To illuminate, the first equality follows from the discussion above; we have blown up $(\PP^{1})^{3}$ along four points, and simply not used $\tilde{p}_{2}$. The second equality follows from Theorem~\ref{thm: p3p1-theorem}. The third equality again follows from the discussion above. The final equality holds as the invariant $\langle p^{6} \rangle_{0,3h}^{\PP^{3}}$ counts the number of degree-3 rational curves in $\PP^{3}$ through six general points. There is only one such curve, the rational normal curve.
\end{example}

\begin{remark}\textnormal{{\bf (Blowing up additional points)}
In the statement of Theorem~\ref{thm: p3p1-theorem}, the the number of points of $\PP^3$ that are blown up is restricted to $6$ to not overburden the notation with indices. In general one may blowup $\PP^3$ at $k\geq 6$ points and $(\PP^1)^3$ at $k-2$ points. The statement of Theorem~\ref{thm: p3p1-theorem} remains true after translating the statement \textit{mutatis mutandis}, setting $\tilde a_i = a_{i+2}$ for $i>3$.}
\end{remark}

\subsection{Strategy of proof} In general, Gromov--Witten invariants are not preserved by birational maps and their transformation under such maps is subtle. Instead, our approach is to relate both geometries to that of the  permutohedral toric threefold. This threefold, denoted $X_{\Pi_3}$, is the blowup of $\PP^3$ at its four coordinate points and six coordinate lines. We prove the equivalence of the all-genus virtual dimension zero non-exceptional GW theories of four spaces, illustrated below.

\begin{figure}[h]
\begin{tikzpicture}[description/.style={fill=white,inner sep=2pt}]
    \matrix (m) [matrix of math nodes, row sep=3em,
    column sep=2.5em, text height=1.5ex, text depth=0.25ex] 
{ GW(\hat{X}) & &  GW(\hat{\tilde{X}}) \\
        GW(X) & & GW(\tilde{X})\\ };
    \path[-,font=\scriptsize]
    (m-1-1) edge node[auto] {Isomorphism} (m-1-3)
            edge node[description] {Blowup} (m-2-1)
    (m-2-1) edge node[auto] {Birational} (m-2-3)
    (m-2-1) edge node[auto,swap] {Transformation} (m-2-3)
    (m-1-3) edge node[description] {Blowup} (m-2-3);
\end{tikzpicture}
\end{figure}

\begin{figure}[h]
\begin{center}
\begin{tikzpicture}[%
  back line/.style={densely dotted}]
  \coordinate (A) at (0,1,1);
  \coordinate (B) at (1,1,1);
  \coordinate (C) at (0,0,1);
  \coordinate (D) at (1,0,1);
 
  \coordinate (A1) at (0,1,0);
  \coordinate (B1) at (1,1,0);
  \coordinate (C1) at (0,0,0);
  \coordinate (D1) at (1,0,0);
 
  \draw[back line] (D1) -- (C1) -- (A1);
  \draw[back line] (C) -- (C1);
  \draw (D1) -- (B1) -- (A1) -- (A)  -- (B) -- (D) -- (C) -- (A);
  \draw (D) -- (D1) -- (B1) -- (B);
  \draw[ultra thick, color=green] (B1)--(B);
  \draw[ultra thick, color=green] (B1)--(A1);
  \draw[ultra thick, color=green] (C1)--(C);
  \draw[ultra thick, color=green] (C)--(A);
  \draw[ultra thick, color=green] (B1)--(D1);
  \draw[ultra thick, color=green] (C)--(D);
  \draw[ball color=blue] (C) circle (.5mm);
  \draw[ball color=blue] (B1) circle (.5mm);
\end{tikzpicture}
\caption{\small{The six lines in $\tilde{X}$} blown up to obtain $\hat{\tilde{X}}$} \label{fig: six-lines-cube}
\end{center}
\end{figure}
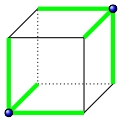

Here $\hat X$ is the blowup of $X$ at six lines and $\hat{\tilde X}$ is the blowup of $\tilde X$ at six lines, see Figure~\ref{fig: six-lines-cube}. Our starting point is the observation that both these varieties come equipped with an explicit isomorphism with the blowup of $X_{\Pi_3}$ at $2$ points. The equivalence of the all genus virtual dimension zero GW theories of $X$ and $\hat{X}$ for nonexceptional (see below) classes was proved by Bryan and the first author~\cite[Lemma 7]{BK}; in this work we complete the square. The non-exceptional property considered here is as follows.

\begin{definition}\label{def: non-exceptional}
Let $\pi :\hat{Y} \rightarrow Y$ be the blowup of the variety $Y$
centered at $Z \subset Y$ with exceptional divisor $E$. A class
 $\beta \in H_{2}(Y;\ZZ )$ is {\emph{nonexceptional}} if:
\begin{enumerate}[(A)]
\item Any stable map to $Y$ representing $\beta$ has image disjoint from the center $Z$;
\item Any stable map to $\hat Y$ representing $\hat{\beta} = \pi^{!}\beta$ has image disjoint the exceptional divisor $E$.
\end{enumerate}
\end{definition}

We prove the equality along the right vertical arrow above by proving that the numerical criterion on $\beta$ appearing in the statement of Theorem~\ref{thm: p3p1-theorem} implies the strong non-exceptional property above. More precisely, we show the following.

\begin{theorem}\label{thm: blowup}
Let $\tilde{X}$ and $\hat{\tilde{X}}$ be as above, and let $\td _{1},\td_{2},\td_{3},\ta_{1},\dotsc ,\ta_{4}\in \ZZ$ be such that $\sum_{}\td_{j}=\sum_{}\ta_{i}$ and $\ta_{i}\neq 0$ for some $i>2$. Then for $\tilde{\beta} = \sum_{}\td_{i}\tilde{h}_{i} - \sum_{}\ta_{j}\tilde{e}_{j}$, we have
\[
\langle \; \rangle^{\tilde{X}}_{g,\tilde{\beta}} =
\langle \; \rangle^{\hat{\tilde{X}}}_{g,\hat{\tilde{\beta}}}, 
\]
where $\hat{\tilde{\beta}} = \sum_{}\td_{i}\hat{\tilde{h}}_{i}-\sum_{}\ta_{j}\hat{\tilde{e}}_{j}$.
\end{theorem}

The symmetry of invariants on $\tilde{X}$ is then deduced by first observing a symmetry on the permutohedral threefold $X_{\Pi_3}$. This induces a symmetry on the invariants of $\hat{\tilde X}$ and hence of $\tilde X$. 

\subsection{Related work and further discussion}\label{sec: further-discussion} View $X_{\Pi_n}$ as the blowup of $\PP^n$ at its coordinate subspaces. The associated fan $\Sigma_{\Pi_n}$ is symmetric about the origin, and the induced automorphism of $X_{\Pi_n}$ resolves the classical Cremona birational map on $\PP^n$. The Cremona symmetry on the Gromov--Witten theory of $\PP^2$ was observed and studied independently in~\cite{CrMi,GoPa}. The Cremona symmetry on $\PP^3$ was studied in the context of the closed topological vertex~\cite{BK}. We remark that such toric symmetries seem to be relatively rare phenomena -- by the results of~\cite{KRRW}, the Cremona symmetry is the only nontrivial toric symmetry that gives rise to relations between the invariants of $\PP^n$ for $n = 2,3$. The analogous symmetry on the genus $0$ Gromov--Witten theory of $\PP^n$ has been studied in~\cite{GKP} using degeneration techniques. Theorem~\ref{thm: p1-symmetry-theorem} uncovers a similar but distinct symmetry on the Gromov--Witten theory of $(\PP^1)^3$. It would be interesting to study to what extent the results of this paper can be generalized to higher genus and higher dimensions. 

\subsection*{Acknowledgements}
D.K. was partially supported by the Beckman Foundation. D.R. was an undergraduate at Harvey Mudd College, supported by a Borrelli Fellowship, when this work began, and it is a pleasure to acknowledge the institution here. We thank the anonymous referee for numerous comments that helped improve the article. 

\section{Gromov-Witten theory}\label{sec: GW} We now briefly recall GW theory and fix notation. Let $X$ be a smooth complex projective variety, and let $\beta \in A_1(X)$ be a curve class. We denote by $\M_{g,n}(X,\beta)$ the moduli stack of stable maps 
\[
f:(C,p_{1},\dotsc ,p_{n}) \to X,
\]
where $C$ is an $n$-marked, possibly nodal genus $g$ curve. This stack admits a virtual fundamental class $[\overline{M}_g(X,\beta)]^{\textnormal{vir}}$ of algebraic degree
\[
\textnormal{vdim}(\M_{g,n}(X,\beta)) = (\textnormal{dim}\ X-3)(1-g)-K_X\cdot \beta +n .
\]
Denote by $ev_{i}$ evaluation morphisms $ev_{i}: \M_{g,n}(X, \beta) \rightarrow X$ defined by 
\[
(f,C,p_{1},\dotsc ,p_{n}) \mapsto f(p_{i}).
\]

Let $\gamma_{1},\dotsc ,\gamma_{n}\in H^{*}(X)$ be a collection of cohomology classes. The genus-$g$ class $\beta $ Gromov-Witten invariant of $X$ with insertions $\gamma_{i}$ is defined by
\[
\langle \gamma_{1},\dotsc ,\gamma_{n} \rangle^X_{g,\beta} = \int_{[\overline{M}_g(X,\beta)]^{\textnormal{vir}}} \prod_{i=1}^{n} ev_{i}^{*}\gamma_{i}.
\]
For further details regarding the fundamentals of Gromov-Witten theory see, for instance~\cite{Hori-et-al}.

\section{Toric Blowups and the Permutohedron}\label{sec: permutohedron}

In this section we construct $\hat{X}$ and $\hat{\tilde{X}}$. First consider the case $k=4$, i.e. we blowup $\PP^{3}$ at four points and $(\PP^{1})^{3}$ at two points. In this case, we are in the toric setting. We prove that $\hat{X}$ and $\hat{\tilde{X}}$ are isomorphic, and in fact are both isomorphic to the permutohedral variety. It follows that $\hat{X } \cong \hat{\tilde{X}}$ for general $k>4$ by simply blowing up along additional points, which need not be fixed.

From the viewpoint of the dual polytopes of these varieties, to construct $\hat{X}$ we realize the permutohedron as a truncation of the simplex, which is classical. However the permutohedron is also constructible by truncation of the cube, yielding $\hat{\tilde{X}}$. This construction is not original; for example Devadoss and Forcey~\cite{Dev08} use this truncation of the cube to construct the permutohedron.

\subsection*{Notation}\label{subsec: notation}
Let $Y$ be a toric variety with fan $\Sigma_Y$. We will denote torus fixed subvarieties in multi-index notation corresponding to generators of their cones. For instance, $p_{i_1\ldots i_k}$ will denote the torus fixed point which is the orbit closure of the cone $\sigma=\langle v_{i_1},\ldots,v_{i_k}\rangle$, for $v_i\in \Sigma_Y^{(1)}$. Similarly $\ell_{i_1\ldots i_r}$ will denote the line which is the orbit closure of $\sigma=\langle v_{i_1},\ldots,v_{i_r}\rangle$, and so on. Further, $Y(Z_1,\ldots,Z_s)$ will denote the iterated blowup of $Y$ at the subvarieties $Z_1,\ldots,Z_s$. By abuse of notation, we will denote by $Y(k)$ the blowup of $Y$ at $k$ points. 

\subsection{The fans of $\hat{X}$ and $\hat{\tilde{X}}$}\label{subsec: fans} The fan $\Sigma_{\PP^3}\subset \ZZ^3$ of $\PP^3$ has 1-skeleton with primitive generators
\[
\left.
\begin{array}{l@{\qquad}c}
v_1 = (-1,-1,-1)  &  v_2 = (1,0,0)\\
 v_3 = (0,1,0) & v_4 = (0,0,1),
\end{array}
\right.
\]
and maximal cones given by 
\[
\left.
\begin{array}{cc}
\langle v_1,v_2,v_3 \rangle  &  \langle v_1,v_2,v_4 \rangle\\
\langle v_1,v_3,v_4 \rangle & \langle v_2,v_3,v_4 \rangle.
\end{array}
\right.
\]
Also note that the fan $\Sigma_{(\PP^1)^{\times3}}\subset \ZZ^3$ of $(\PP^1)^{\times3}$, has primitive generators
\[
\left.
\begin{array}{l@{\qquad}cc}
u_1 = (1,0,0)  &  u_3 = (0,1,0) & u_5 = (0,0,1)\\
u_2 = (-1,0,0) & u_4 = (0,-1,0) & u_6 = (0,0,-1),
\end{array}
\right.
\]
and maximal cones given by 
\[
\left.
\begin{array}{cccc}
\langle u_1,u_3,u_5 \rangle  &  \langle u_1,u_2,u_4 \rangle & \langle u_1,u_2,u_3 \rangle &  \langle u_1,u_2,u_4 \rangle\\
\langle u_2,u_4,u_6 \rangle  & \langle u_2,u_3,u_4 \rangle & \langle u_1,u_2,u_3 \rangle  &  \langle u_1,u_2,u_4 \rangle.
\end{array}
\right.
\]
The three dimensional permutohedron $\Pi_3$ is precisely realized as the dual polytope of the blowup of $\PP^3$ at its 4 torus fixed points and the 6 torus invariant lines between them,
\[
X_{\Pi_3} \cong \PP^3(p_{123},p_{124},p_{134},p_{234},\ell_{12},\ell_{13},\ell_{14},\ell_{23},\ell_{24},\ell_{34}). 
\] 
It is also realized as the dual polytope of a blowup of $(\PP^1)^{\times3}$. In particular,
\[
X_{\Pi_3} \cong (\PP^1)^{\times3}(p_{135},p_{246},\ell_{13},\ell_{15},\ell_{35},\ell_{24},\ell_{26},\ell_{46}).
\]
This blowup of $(\PP^{1})^{\times 3}$ can be viewed as the blowup of two antipodal vertices on the 3-cube and the 6 invariant lines intersecting these points, as shown in Figure~\ref{fig: six-lines-cube}. This common blowup yields an isomorphism $\hat{\tau}:\hat{X}\to \hat{\tilde{X}}$ and a birational map $\tau : \PP^3(4) \rightarrow (\PP^1)^{3}(2)$.
\begin{figure}[h!]
\begin{tikzpicture}[description/.style={fill=white,inner sep=2pt}]
    \matrix (m) [matrix of math nodes, row sep=3em,
    column sep=2.5em, text height=1.5ex, text depth=0.25ex]
    { \hat{X} & & \hat{\tilde{X}} \\
       \PP^3(4) & & (\PP^1)^{3}(2)  \\ };
    \path[->,font=\scriptsize]
    (m-1-1) edge node[auto] {\tiny{$\hat\tau$}} (m-1-3)
    (m-1-3) edge node[description] {$ \tilde{\pi} $} (m-2-3)
    (m-1-1) edge node[description] {$\pi$} (m-2-1);
    \path[->,densely dashed] (m-2-1) edge node[auto] {\tiny{$\tau$}} (m-2-3);
\end{tikzpicture}
\caption{Two constructions of $X_{\Pi_3}$ as a blowup.}
\label{fig: 3-permutohedron-blowups}
\end{figure}


\subsection{Chow Rings}\label{subsec: chow} 
\subsection*{Notation}\label{subsec: more-notation}
We will use $D_{\alpha}$ for the divisor class corresponding to $v_\alpha$ or $u_\alpha$. For blowups, we will label a new element of the 1-skeleton, introduced to subdivide the cone $\sigma = \langle v_i,\ldots,v_j \rangle$, by $v_{i \cdots j}$. Foundations of this material may be found, for instance, in~\cite{Fu}. Classes on $\PP^3(k)$ remain undecorated, tilde classes, such as $\tilde H_i$ or $\tilde e_{ijk}$ signify classes on $(\PP^1)^{3}(k)$, and classes pulled back to the blowup $X_{\Pi_3}$ will be decorated with a hat.


\subsubsection{$X_{\Pi_{3}}$ as a Toric Blowup of $\PP^3$} The Chow ring of $\PP^3$ is generated by the first Chern class of hyperplane bundle  on $\PP^3$. Let $\hat H$ be the pullback of this class to $X_{\Pi_3}$ and let $ \hat{h}= \hat H\cdot \hat H$ denote the class of a general line in $A_1(X)$. Let $\hat E_\alpha$ be the class of the exceptional divisor above the blowup of $p_\alpha$, and $\hat e_\alpha$ be the line class in the exceptional divisor. Let $\hat F_{\alpha'}$ denote the class of the exceptional divisor above the blowup of the line $\ell_{\alpha'}$. Note that that this divisor is abstractly isomorphic to $\PP^1\Cross\PP^1$, so we let $\hat f_\alpha$ and $\hat s_\alpha$ be the section and fiber class respectively. Observe that 
\[
A_2(X_{\Pi_3}) = \langle \hat H, \hat E_{\alpha}, \hat F_{\alpha'} \rangle, \ \ A_1(X_{\Pi_3}) = \langle \hat h, \hat e_{\alpha}, \hat f_{\alpha'} \rangle.
\]
The divisor classes corresponding to $\Sigma_{X_{\Pi_3}}^{(1)}$, are written in terms of this basis as
\begin{eqnarray*}
D_i &=& \hat H - \sum_{i\in \alpha} \hat E_{\alpha} - \sum_{j\in \alpha'}\hat F_{\alpha'} \\
D_{ij} &=& \hat F_{ij}\\
D_{ijk} &=& \hat E_{ijk}. 
\end{eqnarray*}

\subsubsection{As a Toric Blowup of $(\PP^1)^{3}$}\label{subsec: main-map} Let $\hat{\tilde H}_1$, $\hat{\tilde H}_2$ and $\hat{\tilde H}_3$ be the 3 hyperplane classes pulled back from the K\"unneth decomposition of the homology of $(\PP^1)^{3}$. We let $\hat h_{ij}$ be the line class $\hat{\tilde H}_{i}\cdot \hat{\tilde H}_j$ and $\hat{\tilde E}_{\alpha}$, $\hat{\tilde e}_{\alpha}$, $\hat{\tilde F}_{\alpha'}$, $\hat{\tilde f}_{\alpha'}$ and $\hat{\tilde s}_{\alpha'}$ be as above. These classes generate the Chow groups in the appropriate degree. The divisor classes corresponding to $\Sigma_{X_{\Pi_3}}^{(1)}$ are given by
\[
\left.
\begin{array}{cc}
D_1 = \hat {\tilde H}_1 - \hat{\tilde E}_{135} - \hat{\tilde F}_{13}-\hat{\tilde F}_{15},  &  D_2 = \hat {\tilde H}_3 - \hat{\tilde E}_{246} - \hat{\tilde F}_{24}-\hat{\tilde F}_{26}\\

D_3 = \hat {\tilde H}_2 - \hat{\tilde E}_{135} - \hat{\tilde F}_{13}-\hat{\tilde F}_{35}, & D_4 = \hat {\tilde H}_2 - \hat{\tilde E}_{246} - \hat{\tilde F}_{24}-\hat{\tilde F}_{46}\\
D_5 = \hat {\tilde H}_3 - \hat{\tilde E}_{135} - \hat{\tilde F}_{13}-\hat{\tilde F}_{25},  &  D_6 = \hat {\tilde H}_3 - \hat{\tilde E}_{246} - \hat{\tilde F}_{26}-\hat{\tilde F}_{46}\\
D_{ijk} = {\hat{\tilde E}}_{ijk}, & D_{ij} = {\hat {\tilde F}}_{ij}. 
\end{array}
\right.
\]

The map $\hat\tau: \hat{X}\rightarrow \hat{\tilde{X}}$, introduced in Figure~\ref{fig: 3-permutohedron-blowups},  is an isomorphism induced by a relabeling of the fan $\Sigma_{X_{\Pi_3}}$. In particular, the action of $\hat\tau_\star$ on $A_1(X_{\Pi_3})$, is given by 
\begin{eqnarray*}
\hat\tau_\star \hat h &=& \hat {\tilde h}_{12}+ \hat {\tilde h}_{13}+ \hat {\tilde h}_{23}-\hat {\tilde e}_{246}\\
\hat\tau_\star \hat e_{123} &=& \hat {\tilde h}_{13}+ \hat {\tilde h}_{23}-\hat {\tilde e}_{246}\\
\hat\tau_\star \hat e_{124} &=& \hat {\tilde h}_{12}+ \hat {\tilde h}_{23}-\hat {\tilde e}_{246}\\
\hat\tau_\star \hat e_{134} &=& \hat {\tilde h}_{12}+ \hat {\tilde h}_{13}-\hat {\tilde e}_{246}\\
\hat\tau_\star \hat e_{234} &=& {\hat{\tilde e}}_{135}\\
\hat\tau_\star \hat f_{12} &=& \hat {\tilde s}_{46} = \hat {\tilde h}_{23}-\hat{\tilde e}_{246}+\hat{\tilde f}_{46}\\
\hat\tau_\star \hat f_{13} &=& \hat {\tilde s}_{26} = \hat {\tilde h}_{13}-\hat{\tilde e}_{246}+\hat{\tilde f}_{26}\\
\hat\tau_\star \hat f_{14} &=& \hat {\tilde s}_{24} = \hat {\tilde h}_{12}-\hat{\tilde e}_{246}+\hat{\tilde f}_{24}\\
\hat\tau_\star \hat f_{34} &=& {\hat {\tilde f}}_{35}\\
\hat\tau_\star \hat f_{24} &=& {\hat {\tilde f}}_{15}\\
\hat\tau_\star \hat f_{23} &=& {\hat {\tilde f}}_{13}.
\end{eqnarray*}

\section{Toric Symmetries of $\PP^3$ and $(\PP^1)^{3}$}\label{sec: cremona}
The classical Cremona transformation is the rational map 
\[
\xi: \PP^3 \dasharrow \PP^3
\]
defined by 
\[
(x_0:x_1:x_2: x_3) \mapsto (x_1x_2x_3:x_0x_2x_3:x_0x_1x_3:x_0x_1x_2). 
\]
Note that $\xi$ is undefined on the union of the torus invariant points and lines, and is resolved on the maximal blowup of $\PP^3$, $\pi:X_{\Pi_3}\to \PP^3$. The resolved Cremona involution on $X_{\Pi_3}$ is a toric symmetry induced by the reflecting $\Pi_3$ through the origin. Note that the resolved Cremona map, $\hat \xi: \hat{X}\rightarrow \hat{X}$ acts nontrivially on $A^\star(\hat{X})$. For a more detailed treatment of toric symmetries in general and Cremona symmetry in particular, see~\cite{BK, Ga01, KRRW}. Cremona symmetry is given as follows. 
\begin{lemma}[Bryan-Karp~\cite{BK}, Gathmann~\cite{Ga01}]\label{lemma: p3-permutohedron}
Let $\hat{X}$ be the permutohedral blowup of $\PP^3$. Let $\beta$ be  given by
\[
 \beta = d \hat h - \sum_{i=1}^{4} a_i \hat e_i - \sum_{1\leq i< j\leq 6} b_{ij} \hat f_{ij} \in H_{2}(X;\ZZ).
\]
There exists a toric symmetry $\hat \xi$, resolving $\xi$, such that $\hat \xi_\star \beta =\beta '$,
where $\beta' = d'\hat h - \sum_{i} a_i'\hat e_i - \sum_{ij} b_{ij}'\hat f_{ij}$ has coefficients given by 
\begin{eqnarray*}
d' &=& 3d - 2\sum_{i=1}^{4} a_i \\
a_i' &=& d - a_j - a_k - a_l-b_{ij}-b_{ik}-b_{il}\\
b_{ij}' &=& b_{kl},
\end{eqnarray*}
where $\{i,j,k,l\} = \{1,2,3,4\}$.
\end{lemma}
In similar vein, the blowup $\hat{\tilde{X}} \to (\PP^1)^{3}$, also has a nontrivial toric symmetry analogous to Cremona involution. Consider the rational map
\[
\zeta: (\PP^1)^{3} \dasharrow (\PP^1)^{3}
\]
defined by 
\[
((x_0:x_1),(y_0:y_1),(z_0,z_1))\mapsto ((x_1y_0z_0:x_0y_1z_1),(y_0:y_1),(z_0,z_1)).
\]
\begin{lemma}\label{lemma: cube-cremona}
Let $\beta = \sum_{1}^{3} d_{j}\hat{\tilde h}_{j}-a_1\hat{\tilde e}_1-a_2\hat{\tilde e}_2-\sum_{i=1}^6 b_i\hat{ \tilde f}_i\in A_\star (X_{\Pi_3})$. $\hat{\tilde{X}}$ admits a nontrivial toric symmetry $\hat \zeta$, which is a resolution of $\zeta$, whose action on homology is given by 
\[
\hat\zeta_\star \beta = \beta'
\]
where $\beta' = \sum_{1}^{3} d'_{j}\hat{\tilde h}_{j}-a'_1\hat{\tilde e}_1-a'_2\hat{\tilde e}_2-\sum_{i=1}^6 b'_i\hat{ \tilde f}_i$ has coefficients given by
\begin{eqnarray*}
d_{1}' &=& d_{1}+d_{3}-a_1-a_2-b_2-b_5\\
d_{2}' &=& d_{2}+d_{3}-a_1-a_2-b_1-b_4\\
d_{3}' &=& d_{3}\\
a_1'&=& d_{3}-a_2-b_4-b_5\\
a_2' &=& d_{3}-a_1-b_1-b_2\\
b'_1 &=& b_5, \ \ b'_2 = b_4\\
b'_3 &=& b_3, \ \ b'_4 = b_2\\
b'_5 &=& b_1, \ \ b'_6 = b_6.
\end{eqnarray*}
\end{lemma}

\begin{remark}\label{rem: toric-symmetry}
\textnormal{In~\cite{KRRW} it is shown that $X_{\Pi_{3}}$, as a blowup of $\PP^{3}$,  admits a unique nontrivial toric symmetry. Indeed, although the permutohedron admits many symmetries, in the cohomology basis induced by the isomorphism $X_{\Pi_{3}}\cong \hat{X}$, each symmetry is either acts trivially, or is equal to the Cremona symmetry above. Here we have a new toric symmetry of the permutohedron, nontrivial in the cohomology basis induced by $X_{\Pi_{3}}\cong \hat{\tilde{X}}$.}
\end{remark}

\textsc{Proof.}
Observe that choosing $\hat \zeta$ to be the toric symmetry 
\[
\hat\zeta = \begin{bmatrix}
-1 & 0 & 0\\
-1&1&0\\
-1&0&1
\end{bmatrix},
\]
$\zeta_\star$ on $A_\star(X_{\Pi_3})$ has the desired action on homology, and the natural blowup-blowdown composition 
\[
\begin{tikzpicture}[description/.style={fill=white,inner sep=2pt}]
    \matrix (m) [matrix of math nodes, row sep=3em,
    column sep=2.5em, text height=1.5ex, text depth=0.25ex]
    { X_{\Pi_3} & & X_{\Pi_3} \\
       (\PP^1)^{ 3} & & (\PP^1)^{ 3}  \\ };
    \path[->,font=\scriptsize]
    (m-1-1) edge node[auto] {\tiny{$\hat \zeta$}} (m-1-3)
    (m-1-3) edge  (m-2-3)
    (m-1-1) edge  (m-2-1);
    \path[->,densely dashed] (m-2-1) edge node[auto] {\tiny{$\zeta$}} (m-2-3);
\end{tikzpicture}
\]
gives the birational map $\zeta$. 
\qed 


\section{Proof of Main Results}\label{sec: descent}

Recall the notation fixed in the previous section: given a variety $Y$, we denote the blowup of $Y$ at $m$ general points by $Y(m)$. \smallskip

\noindent
\textsc{Reduction to Theorem~\ref{thm: blowup}.} In Section~\ref{sec: permutohedron} we constructed an explicit isomorphism between $\hat{X}$ and $\hat{\tilde X}$. We now recall the following result of Bryan and the first author, which proves the desired relations between the invariants of $\hat X$ and $X$. 

\begin{theorem}[\cite{BK}, Lemma 7]\label{thm: BK}
Let $d, a_{1},\dotsc ,a_{6}\in \ZZ $ be such that $2d= \sum a_{i}$ and $a_{i}\neq 0$ for some $i>4$. Then
\[
\langle \; \rangle^{X}_{g,\beta} = \langle \; \rangle^{\hat{X}}_{g,\hat{\beta }},
\]
where $\beta = dh - \sum_{} a_{i}e_{i}$ and $\hat{\beta} = d\hat{h} - \sum_{ }a_{i}\hat{e}_{i}$. 
\end{theorem}

The GW invariants of $\hat{X}$ and $\hat{\tilde X}$ coincide since these spaces are isomorphic. Thus, Theorem~\ref{thm: p3p1-theorem} will follow immediately from Theorem~\ref{thm: blowup} -- the equality of non-exceptional invariants of $\hat{\tilde X}$ and $\tilde X$. 

We perform a similar reduction on Theorem~\ref{thm: p1-symmetry-theorem}. Recall that we denoted by $\hat f_\alpha$ the fiber classes in the exceptional divisors above lines in $\tilde X$. Observe from Lemma~\ref{lemma: cube-cremona} that under the symmetry $\hat\zeta_\star$, the class $\hat f_\alpha$ is mapped to $\hat f_{\alpha'}$. To prove that $\hat \zeta_\star$ descends to a symmetry on the invariants $\tilde X$, it suffices to prove Theorem~\ref{thm: blowup}. \smallskip

\noindent
\textsc{Proof of Theorem~\ref{thm: blowup}.} We follow the toric construction of Section~\ref{sec: permutohedron}, blowing up $2$ additional points in general position. That is,
\[
\tilde \pi: X_{\Pi_3}(2) = \hat{\tilde X}\to\tilde X =(\PP^1)^{3}(4). 
\]
Let $\hat \beta = \sum d_{j} \hat{\tilde h}_{j} -\sum a_i \hat{\tilde e}_i$ with $a_i\neq 0$ for some $i>2$. We will check the conditions of Definition~\ref{def: non-exceptional} to conclude the theorem. More specifically, we will argue that any stable map in the isomorphism class $[\hat f]\in \overline{M}_g(\hat X, \hat \beta)$ has an image disjoint from $F = \cup \hat {\tilde F}_{jk}$ where the union is taken over all the exceptional divisors above line blowups. We similarly show that any stable map $[f]\in \overline{M}_g(\tilde X,\beta)$ has an image disjoint from $\ell = \cup \ell_{jk}$. It then follows that the map on moduli stacks induced by $\tilde \pi$ is an isomorphism of moduli stacks, compatible with obstruction theories, and hence induces an isomorphism on virtual classes.

Let $[f:C\to \tilde X]\in \overline{M}_g(\tilde X,\beta)$ be a stable map. Suppose that $\textnormal{Im}(f) \cap \ell_{rs} \neq \emptyset$ where $\ell_{rs}$ is one of the six lines in the exceptional locus. Without loss of generality, we assume $\tilde a_3\neq 0$. In particular, the image of $f$ is not completely contained in the line $\ell_{rs}$. Hence we may express the class of the image as 
\[
f_\star [C]= C' + b \ell_{rs}, \ \ (b\geq 0).
\]
Here $C'$ meets $\ell_{rs}$ at finitely many points for topological reasons. Let $\hat C'$ be the $\tilde \pi$-proper transform of $C'$. Since $C'\cap \ell_{rs}\neq \emptyset$, we have $\hat C'\cdot \hat F_{rs} = m>0$. Thus, we write 
\[
\hat C' = \hat \beta - b(\hat{\tilde h}_{j}-\hat{\tilde e}_{\alpha}) - m\hat{\tilde f}_{rs}.
\]
Here $\alpha\in \{1,2\}$, or in other words, $e_{\alpha}$ is the exceptional line above one of the torus fixed points, and $[\ell_{rs}] = \tilde h_{j}$. Apply $\hat{\tau}^{-1}_{\star}$ to the class $\hat C'$ as described explicitly in Section~\ref{subsec: main-map}. We obtain a curve class on $\hat X$:
\[
\hat{\tau}_{\star}^{-1} [\hat C'] = d \hat h - \sum a_i \hat e_i - b(\hat h-\hat e_{\gamma}-\hat e_{\delta}) - m \hat f_{pq},
\]
where $\{\gamma,\delta\}\subset \{1,2,3,4\}$. 
Now consider the divisor 
\[
\hat D_{pq} = 2 \hat H - \sum E_i -\hat F_{pq}-\hat F_{p'q'},
\]
where $\{p,q,p',q'\} = \{1,2,3,4\}$. We calculate that the intersection number of $\hat\tau^{-1}_\star[\hat C]$ with $D_{pq}$ is $-m$, however from~\cite[Lemma 11]{BK}, $\hat D_{pq}$ is nef so we have a contradiction. We conclude that $\textrm{Im}(f) \cap \ell_{rs} = \emptyset$. This proves that the classes $\beta$ as given, satisfy condition (A) of Definition~\ref{def: non-exceptional}. 

We argue in similar fashion for $\overline{M}_g(\hat X,\hat \beta)$. Let $[\hat f: C\to \hat X]$. Suppose $\textnormal{Im}(\hat f)\cap \hat{\tilde F}_{rs}\neq \emptyset$. Since $\hat \beta \cdot \hat{\tilde F}_{rs} = 0$, $f_\star C$ must have a component $C''$ completely contained in $\hat{\tilde F}_{rs}$, where we have 
\[
f_\star C = C'+C'',
\]
where $C'$ is nonempty since $\hat \beta \cdot \hat{\tilde E}_{4}\neq 0$. Since $C''\subset \hat {\tilde F}_{rs}$ is an effective class in $\hat {\tilde F}_{rs} \cong \PP^1\Cross \PP^1$, it must be of the form $C'' = a \hat f_{rs}+b \hat s_{rs}$ for $a,b\geq 0$ and $a+b >0$. We compute  $\hat{\tau}_{\star}^{-1} (\hat D_{pq}) \cdot C' = -a-b$, contradicting the fact that $\hat{D}_{pq}$ is nef. Again we may conclude that $\textnormal{Im}(\hat f)\cap \hat {\tilde F}_{rs}= \emptyset$, and $\beta$ satisfies condition (B) of Definition~\ref{def: non-exceptional}.\qed  

\bibliographystyle{siam}
\bibliography{P1xP1xP1}

\end{document}